\documentclass[english]{article}
\usepackage[T1]{fontenc}
\usepackage[latin9]{inputenc}
\usepackage{esint}

\makeatletter
\newcommand{\lyxaddress}[1]{
\par {\raggedright #1
\vspace{1.4em}
\noindent\par}
}

\makeatother

\usepackage{babel}
\begin{document}

\title{Multitimescale method for approximating the path action relevant
to non-equilibrium statistical physics.}

\author{Richard Kleeman}
\maketitle

\lyxaddress{Courant Institute of Mathematical Sciences, New York.}
\begin{abstract}
A path integral formalism has been proposed recently for non-equilibrium
statistical physics applications (see \cite{Kle14}). In this contribution
we outline an efficient method for its numerical evaluation. The method
used is based on the multiscale MCMC method of Ceperley and co-workers
in quantum applications. A significant new feature of the method proposed
is that the time endpoint is not fixed and indeed the endpoint sample
is the principal object of interest.
\end{abstract}

\section{Background}

The action proposed in \cite{Kle14} may be written as
\begin{equation}
S=\Delta t\int\left\{ \left(\dot{\lambda}-\Theta\left(\lambda\right)\right)^{t}g\left(\dot{\lambda}-\Theta\left(\lambda\right)\right)+\Psi\left(\lambda\right)\right\} dt\label{contaction}
\end{equation}

where $\Psi\geq0$ and $g$ is a Riemannian (non-negative definite)
tensor often called the Fisher information matrix. $\Delta t$ is
the minimum timescale associated with the slow modes of the system
being modelled. A generalized Boltzmann principle then constructs
a path measure using $\exp\left(-S\right)$. Quantities of interest
for this theory can be obtained by constructing a sample set of paths
according to the proposed measure. Such a sample can be obtained using
the Markov Chain Monte Carlo methodology pioneered in the context
of quantum Bose condensates by Ceperley and co-workers \cite{ceperley1995path}.
This uses a multitimescale approach to speed path sampling via a Brownian
bridge style interleaving of time nodes. Thus at the coarsest level
the midpoint between the assumed fixed endpoints is retained. At the
next level the midpoints of all three previously defined nodes are
used. This ``interleaving'' of nodes continues until the finest
level desired is obtained. Ceperley proposes that at each multiscale
level (bar the finest) an approximate rapidly computable action $S_{k}=-\log\pi_{k}$
be used. This is a function of the variables at level $k$ and below
which are denoted by $s_{0},s_{1},\ldots,s_{k}$ and abbreviated as
$s$. 

As is usual in MCMC methods, in order to make a transition at level
$k$ of $s_{k}\rightarrow s'_{k}$ we require a trial transition probability
function at each level for which samples are rapidly obtainable. Note
that transitions are made starting at the coarsest level and proceeding
to the finest. This is facilitated by the interleaving of multiscale
time points. Denote this transition probability for level $k$ by
$T_{k}(s'_{k})$ and note that it will not, in our case, depend on
the starting state $s_{k}$. Ceperley then shows that if we choose
the acceptance rule
\begin{equation}
A_{k}(s_{k}')=\min\left\{ 1,\frac{T_{k}(s_{k})\pi_{k}(s')\pi_{k-1}(s)}{T_{k}(s_{k}')\pi_{k}(s)\pi_{k-1}(s')}\right\} \label{acceptance}
\end{equation}

as well as ensuring that the finest $\pi_{m}$ is exact (up to time
discretization) then a sample will be produced according to $\pi_{m}$
and the multiscaling will aid in producing a ``reasonably rapid''
sampling of path space. 

Let us suppose we were able to select $\pi_{k}$ to be the marginal
density with respect to the variables $s_{0},s_{1},\ldots,s_{k}$
i.e. the variables $s_{k+1},\ldots,s_{m}$ are integrated out from
$\pi_{m}$. Then if we select
\[
T_{k}^{*}(s_{k})=\pi_{k}/\pi_{k-1}
\]

it is easily seen firstly that $T_{k}^{*}$ is actually a probability
density for the $s_{k}$ (since integrating out the $s_{k}$ of $\pi_{k}$
obviously produces $\pi_{k-1}$ so it is normalized correctly) and
secondly with such a choice all trial transitions are accepted (see
equation (\ref{acceptance})) at all levels meaning a very rapid algorithm
since in MCMC rejected transitions imply further sampling of the original
$s$.

The above strategy is called the heat bath MCMC method. Unfortunately,
of course, it is only practical for very special $\pi_{m}$ i.e. ones
which may be integrated analytically. The best known of these are
the Gaussian densities. For the application under consideration the
actions are never exactly Gaussian so instead we produce a good Gaussian
approximation $\pi_{m}^{a}$ to the exact $\pi_{m}$ and use this
to produce, by analytical integration, a set of $\pi_{k}^{a}$. We
then set
\begin{eqnarray*}
\pi_{k} & = & \pi_{k}^{a}\quad k<m\\
T_{k} & = & \pi_{k}^{a}/\pi_{k-1}^{a}\quad k\leq m
\end{eqnarray*}

With such a choice it is clear that acceptance always occurs at all
levels bar the finest when the exact $\pi_{m}$ is required in the
acceptance rule (\ref{acceptance}). For this level we obtain
\[
A_{m}(s')=\min\left(1,\frac{\pi_{m}^{a}(s)\pi_{m}(s')}{\pi_{m}^{a}(s')\pi_{m}(s)}\right)
\]

This strategy is quite different to that used in Bose condensates
by Ceperley. There the molecular interaction term in the action is
highly non-quadratic for close approaches of the often densely packed
molecules. In our case a quadratic approximation is responsible for
a large part of the irreversible relaxation of the statistical system
at least for the examples studied to date (see \cite{kleeman2012nonequilibrium}).
The current strategy also has the practical advantage that Gaussian
densities for $T_{k}$ enable rapid sampling. It is reasonably clear
also that the amount of acceptance depends on the accuracy of the
quadratic approximation to the last level action.

\section{Gaussian approximation of action}

There are evidently many ways in which a quadratic approximation of
the action could be constructed. The optimal strategy would be to
ensure that the approximation was most accurate for the most likely
paths. Intuitively this might be ensured if linearizations were performed
about a path close to the mean trajectory of the statistical system.
Since this trajectory is not known a priori, then a judicious guess
would seem the best approach. Note that theoretically the accuracy
of the approximation only affects the speed of convergence of the
MCMC produced sample and not the ultimate accuracy of a reasonably
large sample. Thus if convergence (acceptance) is too slow then this
may be an indication of a suboptimal action approximation. Let us
assume that a trial trajectory $\overline{\lambda}(t)$ is specified
(see below for possible choices). 

We then choose initially to approximate $g$ and $\Psi$ as respectively
$g(\lambda=0)$ and $\phi_{ij}\lambda_{i}\lambda_{j}$. The second
on the basis that the irreversible term is minimized at $\lambda=0$
which from numerical studies is often the equilibrium trial density.
We linearize $\Theta$ about $\overline{\lambda}(t)$ and write
\begin{equation}
\Theta\simeq\overline{A}(t)+A(t)(\lambda-\overline{\lambda})\label{linearize}
\end{equation}

where $\overline{A}$ and $A$ are computible from $\Theta$ and $\overline{\lambda}(t)$.
The time discretized approximate action at level $m$ can then be
written as
\begin{eqnarray}
S(\lambda) & \simeq & \tau_{s}\Delta t\sum_{n=1}^{N}\left(\frac{\lambda_{i}^{n+1}-\lambda_{i}^{n}}{\Delta t}-\frac{1}{2}\left(A_{ik}^{n+1}\lambda_{k}^{n+1}+A_{ik}^{n}\lambda_{k}^{n}\right)-\overline{B}_{i}^{n+0.5}\right)g_{ij}\left(i\leftrightarrow j\right)\nonumber \\
 &  & \qquad+\frac{1}{2}\phi_{ij}\left(\lambda_{i}^{n}\lambda_{j}^{n}+\lambda_{i}^{n+1}\lambda_{j}^{n+1}\right)\label{Gaussapprox}\\
\overline{B}_{i}^{n+0.5} & = & \overline{A}_{i}^{n+0.5}-\frac{1}{2}\left(A_{ik}^{n+1}\overline{\lambda}_{k}^{n+1}+A_{ik}^{n}\overline{\lambda}_{k}^{n}\right)
\end{eqnarray}

where the summation of repeated indices convention is being assumed;
the second bracket is simply the first with $j$ replacing $i$; $\tau_{s}$
is notionally equal to the fine time step $\Delta t$ but could be
subjected to tuning given the uncertainty over the precise value of
the latter quantity. Finally $\overline{A}$ is evaluated at the midpoint
of $n$ and $n+1$. In order to perform Gaussian marginalizing integrations
of various $\lambda^{k}$ we need to rewrite $S$ in terms of quadratic
and linear combinations of such variables. 

To facilitate notational efficiency in manipulating (\ref{Gaussapprox}),
replace the upper index $n$ with $0$, $n+1$ with $+$ and $n-1$
with $-$. Two summands will produce terms of the form $\lambda_{i}^{0}\lambda_{j}^{0}$;
$\lambda_{i}^{0}\lambda_{j}^{+}$; $\lambda_{i}^{0}\lambda_{j}^{-}$
and $\lambda_{i}^{0}$ which are needed to evaluate an integration
of the variables $\lambda_{i}^{0}$ . For the first type, from the
summand involving $n$ and $n+1$ we obtain the terms
\[
\tau_{s}\lambda_{i}^{0}\lambda_{j}^{0}\left[\frac{g_{ij}}{\Delta t}+A_{kj}^{0}g_{ki}+\frac{\Delta t}{4}A_{ki}^{0}A_{lj}^{0}g_{kl}+\frac{\Delta t}{2}\phi_{ij}\right]
\]

From the summand involving $n$ and $n-1$ terms of this form will
also occur but in that case the second term in the square brackets
has a reversed sign so we obtain in total 
\begin{equation}
\tau_{s}\lambda_{i}^{n}\lambda_{j}^{n}\left[\frac{2g_{ij}}{\Delta t}+\frac{\Delta t}{2}A_{ki}^{0}A_{lj}^{0}g_{kl}+\Delta t\phi_{ij}\right]\label{Ginit}
\end{equation}

The second type of (quadratic) term comes only from the summand involving
$n$ and $n+1$:
\begin{equation}
\tau_{s}\lambda_{i}^{0}\lambda_{j}^{+}\left[-\frac{2g_{ij}}{\Delta t}+A_{kj}^{+}g_{ik}-A_{ki}^{0}g_{jk}+\frac{\Delta t}{2}A_{ki}^{0}A_{lj}^{+}g_{lk}\right]\label{cross+init}
\end{equation}

The third type of term comes analogously to the last but from the
summand involving $n$ and $n-1$:
\begin{equation}
\tau_{s}\lambda_{i}^{-}\lambda_{j}^{0}\left[-\frac{2g_{ij}}{\Delta t}+A_{kj}^{0}g_{ik}-A_{ki}^{-}g_{jk}+\frac{\Delta t}{2}A_{ki}^{-}A_{lj}^{0}g_{lk}\right]\label{cross-init}
\end{equation}

Turning now to the linear term we note that they occur due to the
presence of $\overline{B}$ in the discetized action. We keep track
of both the $\lambda^{0}$ and $\lambda^{+}$ terms since the latter
in the summand $n-1$ and $n$ will contribute to the total $\lambda^{0}$
piece. The cross terms give 
\[
\tau_{s}\left[2\left(\lambda_{i}^{0}-\lambda_{i}^{+}\right)g_{ij}\overline{B}_{j}^{0+}+\Delta t\overline{B}_{i}^{0+}g_{ij}\left(A_{jk}^{0}\lambda_{k}^{0}+A_{jk}^{+}\lambda_{k}^{+}\right)\right]
\]

($\overline{B}^{0+}\equiv\overline{B}^{n+0.5}$) which rearranges
to
\[
\tau_{s}\left[\lambda_{i}^{0}\left(2g_{ij}\overline{B}_{j}^{0+}+\Delta tA_{ji}^{0}g_{kj}\overline{B}_{k}^{0+}\right)-\lambda_{i}^{+}\left(2g_{ij}\overline{B}_{j}^{0+}-\Delta tA_{ji}^{+}g_{kj}\overline{B}_{k}^{0+}\right)\right]
\]

Thus the total contribution from the action to the linear $\lambda^{0}$
term is
\[
\tau_{s}\lambda_{i}^{0}\left[2g_{ij}\left(\overline{B}_{j}^{0+}-\overline{B}_{j}^{-0}\right)+\Delta tA_{ji}^{0}g_{kj}\left(\overline{B}_{k}^{0+}+\overline{B}_{k}^{-0}\right)\right]
\]

or
\begin{equation}
2\tau_{s}\lambda_{i}^{0}\left[g_{ij}\left(\overline{B}_{j}^{0+}-\overline{B}_{j}^{-0}\right)+\Delta tA_{ji}^{0}g_{kj}\overline{B}_{k}^{0}\right]\label{Kinit}
\end{equation}

\section{Determination of general level action}

As integration of the variables associated with each level $k$ occurs,
a new marginalized action appears which is used to determine the transition
probability $T_{k-1}$ for the next coarser level. A matrix recursion
scheme is required to determine the Gaussian densities $\pi_{k}$.
Because of the Brownian bridge construction we may write the part
of the level $k$ action involving the $2^{k-1}$ variables to be
integrated as 
\begin{eqnarray}
S(k) & = & \sum_{n=1}^{N}\frac{1}{2}\lambda_{i}^{l}\lambda_{j}^{l}G_{ij}^{l}(k)+\lambda_{i}^{l}\left(\lambda_{j}^{l+\Delta}H_{ij}^{+l}(k)+\lambda_{j}^{l-\Delta}H_{ij}^{-l}(k)+K_{i}^{l}(k)\right)\label{Marginal}\\
l(n,k) & \equiv & 1+\Delta(2n-1)\nonumber \\
N & \equiv & 2^{k-1}\nonumber \\
\Delta & \equiv & 2^{m-k}\nonumber 
\end{eqnarray}

where the $G$, $H$ and $K$ need to be determined recursively. This
structure reflects the interleaving nature of the multiscale construction.
There are (excluding endpoints see below) $N-1$ variables not integrated
at this level. They occur at time indices
\[
l(r,k)+\Delta\qquad r=1,\ldots,N-1
\]

Each acquires matrix contributions from integrated variables at $l(r,k)$
and $l(r+1,k)=l(r,k)+2\Delta$. Using the usual multivariate Gaussian
integration formula we can integrate the first of these variables
(functional form of the density being the Gaussian $\exp(-S(k))$
obtaining the contributions
\begin{equation}
-\frac{1}{2}\left[\lambda_{j}^{l+\Delta}H_{ij}^{+l}(k)+\lambda_{j}^{l-\Delta}H_{ij}^{-l}(k)+K_{i}^{l}(k)\right]\left(G^{l}(k)\right)_{ik}^{-1}\left[i\leftrightarrow k\right]\label{integrmult}
\end{equation}

Another analogous contribution comes from the second integration variable
mentioned above. Comparing this with the form of $S(k-1)$ we read
off the recursion relations for the matrices $G$:
\begin{equation}
G^{l+\Delta}(k-1)=G^{l+\Delta}(k)-H^{+l}(k)^{t}G^{l}(k)^{-1}H^{+l}(k)-H^{-(l+2\Delta)}(k)^{t}G^{l+2\Delta}(k)^{-1}H^{-(l+2\Delta)}(k)\label{Grec}
\end{equation}

Now the only terms of the form $\lambda^{l+\Delta}\lambda^{l-\Delta}$
come from (\ref{integrmult}). Furthermore at level $k-1$ only one
of the variables with time index $l+\Delta$ and $l-\Delta$ will
be an integration variable due to the interleaving construction. This
term therefore is assigned entirely to the time index to be integrated
and by convention 
\begin{equation}
H^{+(l-\Delta)}(k-1)=H^{-(l+\Delta)}(k-1)\label{convention}
\end{equation}

Comparison of this term's form with the form of $S(k-1)$ shows that
\begin{equation}
H^{+(l-\Delta)}(k-1)=-\left(H^{-l}(k)\right)^{t}\left[G^{l}(k)\right]^{-1}H^{+l}(k)\label{Hrec}
\end{equation}

Note that the required $H^{-}$ terms are obtainable from (\ref{convention}).
Finally the linear term involving $K$ is easily seen to imply for
the vectors $K$ 
\begin{equation}
K^{l+\Delta}(k-1)=K^{l+\Delta}(k)-H^{+l}(k)^{t}\left[G^{l}(k)\right]^{-1}K^{l}(k)-H^{-(l+2\Delta)}(k)^{t}\left[G^{l+2\Delta}(k)\right]^{-1}K^{l+2\Delta}(k)\label{Krec}
\end{equation}

Equations (\ref{Grec}), (\ref{Hrec}) and (\ref{Krec}) define the
recursion which allows $T_{k}$ to be constructed since this is simply
the multivariate Gaussian $C\exp\left(-S(k)\right)$. The start of
this recursive chain is obtained by setting $k=m$ and comparing with
the equations (\ref{Ginit}), (\ref{cross+init}), (\ref{cross-init})
and (\ref{Kinit}) from the last section:
\begin{eqnarray*}
G_{ij}^{l}(m) & = & 2\tau_{s}\left[\frac{2g_{ij}}{\Delta t}+\frac{\Delta t}{2}A_{ki}^{l}A_{lj}^{l}g_{kl}+\Delta t\phi_{ij}\right]\\
H_{ij}^{+l}(m) & = & \tau_{s}\left[-\frac{2g_{ij}}{\Delta t}+A_{kj}^{l+1}g_{ik}-A_{ki}^{l}g_{jk}+\frac{\Delta t}{2}A_{ki}^{l}A_{lj}^{l+1}g_{lk}\right]\\
H_{ij}^{-l}(m) & = & \tau_{s}\left[-\frac{2g_{ij}}{\Delta t}+A_{kj}^{l}g_{ik}-A_{ki}^{l-1}g_{jk}+\frac{\Delta t}{2}A_{ki}^{l-1}A_{lj}^{l}g_{lk}\right]\\
K_{i}^{l}(m) & = & 2\tau_{s}\left[g_{ij}\left(\overline{B}_{j}^{l+1/2}-\overline{B}_{j}^{l-1/2}\right)+\Delta tA_{ji}^{l}g_{kj}\overline{B}_{k}^{l}\right]
\end{eqnarray*}

Note that for $k=m$ we have $\Delta=1$.

\section{Temporal boundary points}

For the application here, the starting point is fixed but the endpoint
is not. We deal with that in the context of the interleaving ``Brownian
bridge'' construction by placing the floating endpoint in the level
$0$ set of variables and adding the following contribution to the
approximate action at every level:
\begin{eqnarray*}
S_{end}(k) & = & \frac{1}{2}\lambda_{i}^{M}\lambda_{j}^{M}G_{ij}^{M}(k)+\lambda_{i}^{M}K_{i}^{M}(k)\\
M & \equiv & 2^{m}
\end{eqnarray*}

Note the absence of the cross terms which is caused by truncating
the action at the endpoint which omits the forward term and the inclusion
already of the backward term in the final interior terms. This modification
sets recursion relations for $G^{M}$ and $K^{M}$ but note that the
cross term coefficients $H$ are not affected since the ones required
are already computed in the previous section. We can read off the
required relations from the interior relations (\ref{Grec}) and (\ref{Krec})
in the previous section by dropping the forward terms:
\begin{eqnarray*}
G^{M}(k-1) & = & G^{M}(k)-H^{+(M-\Delta)}(k)^{t}G^{M-\Delta}(k)^{-1}H^{+(M-\Delta)}(k)\\
K^{M}(k-1) & = & K^{M}(k)-H^{+(M-\Delta)}(k)^{t}\left[G^{M-\Delta}(k)\right]^{-1}K^{M-\Delta}(k)
\end{eqnarray*}

In terms of sampling, the starting point (level zero variables) now
occurs with the endpoint random variable which follows a Gaussian
described by $G^{M}(0)$ and $K^{M}(0)$. Intuitively one expects
the mean of this density to vary from close to the initial conditions
for short time actions through to values close to that for the mean
of the equilibrium consistency distribution when the action time interval
is long. Once the endpoint is determined the other random variables
are sampled following the interleaving pattern described in the previous
section. The interleaving nature ensures that all sampling distributions
are independent with respect to their time. If the variables are vectors
however this may not be true with respect to the vector index.

\section{Linearization trajectory}

As mentioned previously $\overline{\lambda}(t)$ is required to compute
the various linearized terms associated with $A$. Experience with
simple turbulence systems suggests the following form as a rough approximation
for the evolution of the slow variable co-ordinates of the system:
\begin{equation}
\dot{\overline{\lambda}}-\Theta\left(\overline{\lambda}\right)=-\frac{1}{2}\nabla\Phi(\overline{\lambda})\label{meanpath}
\end{equation}

To reiterate this choice is not required to be exact, simply to ensure
that the linearization is sufficiently accurate for ``most'' paths
in order to ensure that the acceptance rate of the MCMC method is
reasonable. Clearly the ``acid test'' here is this latter rate and
if it is too low for the choice made in (\ref{meanpath}) then another
choice may be indicated.

\section{A sample linearization: The truncated Burgers model.}

Here we have for various slow variable choices, the complex relation
\begin{eqnarray}
\Theta_{2k-1}(\lambda)+i\Theta_{2k}(\lambda) & = & -\frac{ik}{2}\sum_{k=k_{1}+k_{2}}z_{k_{1}}z_{k_{2}}\label{basictbh}\\
z_{j} & \equiv & \lambda_{2j-1}+i\lambda_{2j}\quad0<j<m\nonumber \\
z_{j}^{*} & = & z_{-j}\nonumber \\
z_{0} & = & 0\nonumber \\
z_{l} & = & 0\quad\left|l\right|>m\nonumber 
\end{eqnarray}

The index $k$ runs $1,\ldots,2m$; while $j$, $k_{1}$ and $k_{2}$
run from $-m,\ldots,m$. Linearizing as 
\[
\lambda_{l}=\overline{\lambda}_{l}+\lambda'_{l}
\]

we obtain
\begin{eqnarray*}
\Theta'_{2k-1}+i\Theta'_{2k}=C_{k} & = & -ik\sum_{k=k_{1}+k_{2}}\overline{z}_{k_{1}}z_{k_{2}}=-ik\sum_{k_{2}=-m}^{m}\overline{z}_{k-k_{2}}z'_{k_{2}}\\
 & \equiv & Z_{kk_{2}}z'_{k_{2}}\\
Z_{kl} & = & -ik\overline{z}_{k-l}
\end{eqnarray*}

where the obvious summation convention is assumed on the second line.
Set the following:
\begin{eqnarray*}
Z_{kl} & = & X_{kl}+iY_{kl}\\
z'_{l} & = & b_{l}+ic_{l}
\end{eqnarray*}

so
\begin{eqnarray*}
b_{-\left|j\right|} & = & b_{\left|j\right|}\\
c_{-\left|j\right|} & = & -c_{-\left|j\right|}
\end{eqnarray*}

then it is easily checked that
\begin{eqnarray*}
\Theta'_{2k-1} & = & X_{kk_{2}}b_{k_{2}}-Y_{kk_{2}}c_{k_{2}}=\sum_{k_{2}=1}^{m}\left(E_{kk_{2}}b_{k_{2}}+F_{kk_{2}}c_{k_{2}}\right)\\
\Theta'_{2k} & = & Y_{kk_{2}}b_{k_{2}}+X_{kk_{2}}c_{k_{2}}=\sum_{k_{2}=1}^{m}\left(E'_{kk_{2}}b_{k_{2}}+F'_{kk_{2}}c_{k_{2}}\right)
\end{eqnarray*}

where
\begin{eqnarray*}
E_{kk_{2}} & = & X_{kk_{2}}+X_{k(-k_{2})}\\
E'_{kk_{2}} & = & Y_{kk_{2}}+Y_{k(-k_{2})}\\
F_{kk_{2}} & = & Y_{k(-k_{2})}-Y_{kk_{2}}\\
F'_{kk_{2}} & = & X_{kk_{2}}-X_{k(-k_{2})}
\end{eqnarray*}

Now since $\lambda_{2j-1}=b_{j}$ and $\lambda_{2j}=c_{j}$ we can
construct the matrix $A$ for equation (\ref{linearize}) i.e. we
have
\[
\Theta'=A\lambda'
\]

with
\begin{eqnarray*}
A_{2r(2s-1)} & = & E'_{rs}\\
A_{2r(2s)} & = & F'_{rs}\\
A_{(2r-1)(2s-1)} & = & E_{rs}\\
A_{(2r-1)(2s)} & = & F_{rs}
\end{eqnarray*}

where $0<r<2m$ but $0<s<m$. Note that the vector $\Theta'$ is double
the length of $\lambda'$ because the sum in (\ref{basictbh}) allows
values of $k$ beyond the possible values for $k_{1}$ and $k_{2}$.

\section{The ``free particle'' case}

Set $A=\phi=0$ and $g=I$. We have then
\begin{eqnarray*}
G^{j}(m) & = & 2\alpha I\\
H^{\pm l}(m) & = & -\alpha I\\
\alpha & \equiv & \frac{2\tau_{s}}{\Delta t}
\end{eqnarray*}

and so
\begin{eqnarray*}
G^{l+1}(m-1) & = & G^{l+1}(m)-H^{+l}(m)^{t}\left[G^{l}(m)\right]^{-1}H^{+1}(m)-H^{-(l+2)}(m)^{t}\left[G^{l+2}(m)\right]^{-1}H^{-(1+2)}(m)\\
 & = & 2\alpha I-2\alpha^{2}\frac{1}{2\alpha}I=\alpha I
\end{eqnarray*}

and
\[
H^{\pm(l\pm1)}(m-1)=-H^{-l}(m)^{t}\left[G^{l}(m)\right]^{-1}H^{+l}(m)=-\alpha^{2}\frac{1}{2\alpha}I=-\frac{\alpha}{2}I
\]

Iterating we obtain easily
\begin{eqnarray*}
G^{j}(k) & = & \frac{\alpha}{2^{k-m-1}}I\\
H^{\pm j}(k) & = & -\frac{1}{2}G^{j}(k)
\end{eqnarray*}

This form for $T$ is discussed by Ceperley \cite{ceperley1995path}.
Effectively it implies increasing the time increment $\Delta t$ by
a factor of two as the timescale level becomes coarser by one. 

\section{Most accurate quadratic approximation for the action.}

In order to increase acceptances in the MCMC method and hence reduce
the computational burden, the best quadratic approximation is needed.
The algorithm detailed above produces a 100\% acceptance rate when
the exact action is quadratic since it is a heat bath method. In order
to produce a quadratic approximation we linearized $\Theta\left(\lambda\right)$.
In the present section we instead construct a second order Taylor
expansion for the action using a mean trajectory of the type discussed
above.

The most general discretised action may be written as
\begin{eqnarray*}
S & = & \Delta t\sum_{n=1}^{N}L(n)dt\\
L(n) & = & \left(\frac{\lambda_{i}(n)-\lambda_{i}(n-1)}{\Delta t}-\frac{1}{2}\left(\Theta_{i}(n)+\Theta_{i}(n-1)\right)\right)h^{ij}(n)\left(j\leftrightarrow i\right)+\Psi(n)\\
h^{ij}(n) & = & \frac{1}{2}\left(g^{ij}(n)+g^{ij}(n-1)\right)
\end{eqnarray*}

where the obvious summation convention holds and $n$ is the time
index for which the convention does not hold in what follows. A second
order Taylor expansion about the mean trajectory $\overline{\lambda}$
gives
\begin{eqnarray*}
S^{a} & = & S\biggr\vert_{\lambda=\overline{\lambda}}+\sum_{n=1}^{N}\frac{\partial S}{\partial\lambda_{k}(n)}\biggr\vert_{\lambda=\overline{\lambda}}\lambda'_{k}(n)+\frac{1}{2}\sum_{n=1}^{N}\sum_{m=1}^{N}\frac{\partial^{2}S}{\partial\lambda_{k}(n)\partial\lambda_{l}(m)}\biggr\vert_{\lambda=\overline{\lambda}}\lambda'_{k}(n)\lambda'_{l}(m)\\
\lambda' & = & \lambda-\overline{\lambda}
\end{eqnarray*}

To facilitate computation we introduce the following symbols
\begin{eqnarray}
T_{i}(n) & \equiv & \frac{\lambda_{i}(n)-\lambda_{i}(n-1)}{\Delta t}-\frac{1}{2}\left(\Theta_{i}(n)+\Theta_{i}(n-1)\right)\nonumber \\
G_{ij}^{+}(n) & \equiv & \frac{2\delta_{ij}}{\Delta t}-\frac{\partial\Theta_{i}(n)}{\partial\lambda_{j}(n)}=2\frac{\partial T_{i}(n)}{\partial\lambda_{j}(n)}\nonumber \\
G_{ij}^{-}(n) & \equiv & \frac{-2\delta_{ij}}{\Delta t}-\frac{\partial\Theta_{i}(n)}{\partial\lambda_{j}(n)}=2\frac{\partial T_{i}(n+1)}{\partial\lambda_{j}(n)}\nonumber \\
R_{i}^{kl}(n) & \equiv & \frac{\partial g^{kl}(n)}{\partial\lambda_{i}(n)}=2\frac{\partial h^{kl}(n)}{\partial\lambda_{i}(n)}=2\frac{\partial h^{kl}(n+1)}{\partial\lambda_{i}(n)}\nonumber \\
V_{kij}(n) & \equiv & \frac{\partial^{2}\Theta_{k}(n)}{\partial\lambda_{i}(n)\partial\lambda_{j}(n)}=-\frac{\partial G_{ki}^{+}(n)}{\partial\lambda_{j}(n)}=-\frac{\partial G_{ki}^{-}(n)}{\partial\lambda_{j}(n)}\nonumber \\
U_{ij}^{kl}(n) & \equiv & \frac{\partial^{2}g^{kl}(n)}{\partial\lambda_{i}(n)\partial\lambda_{j}(n)}=\frac{\partial R_{i}^{kl}(n)}{\partial\lambda_{j}(n)}\label{tensors}
\end{eqnarray}

in which notation we have
\[
S=\Delta tdt\sum_{n=1}^{N}\left\{ T_{i}(n)h^{ij}T_{j}(n)+\Psi(n)\right\} 
\]

Differentiating $S$ once using (\ref{tensors}) and also using the
fact that the metric tensor $g$ is symmetric we obtain
\begin{eqnarray*}
\frac{1}{\alpha}\frac{\partial S}{\partial\lambda_{i}(n)} & = & G_{ki}^{+}(n)h^{kl}(n)T_{l}(n)+\frac{1}{2}R_{i}^{kl}(n)T_{k}(n)T_{l}(n)\\
 &  & +G_{ki}^{-}(n)h^{kl}(n+1)T_{l}(n+1)+\frac{1}{2}R_{i}^{kl}(n)T_{k}(n+1)T_{l}(n+1)\\
 &  & +\frac{\partial\Psi(n)}{\partial\lambda_{i}(n)}\\
\alpha & \equiv & \Delta tdt
\end{eqnarray*}

Differentiating again using (\ref{tensors}) and again using the metric
tensor symmetry we obtain the Hessian matrix
\begin{eqnarray*}
\frac{1}{\alpha}\frac{\partial^{2}S}{\partial\lambda_{i}(n)\partial\lambda_{j}(r)} & = & \delta_{rn}\biggl[-V_{kij}(n)h^{kl}(n)T_{l}(n)+\frac{1}{2}G_{ki}^{+}(n)R_{j}^{kl}(n)T_{l}(n)\\
 &  & +\frac{1}{2}G_{ki}^{+}(n)h^{kl}(n)G_{lj}^{+}(n)+\frac{1}{2}U_{ij}^{kl}(n)T_{k}(n)T_{l}(n)+\frac{1}{2}R_{i}^{kl}(n)G_{kj}^{+}(n)T_{l}(n)\\
 &  & -V_{kij}(n)h^{kl}(n+1)T_{l}(n+1)+\frac{1}{2}G_{ki}^{-}(n)R_{j}^{kl}(n)T_{l}(n+1)\\
 &  & +\frac{1}{2}G_{ki}^{-}(n)h^{kl}(n+1)G_{lj}^{-}(n)+\frac{1}{2}U_{ij}^{kl}(n)T_{k}(n+1)T_{l}(n+1)\\
 &  & +\frac{1}{2}R_{i}^{kl}(n)G_{kj}^{-}(n)T_{l}(n+1)+\frac{\partial^{2}\Psi(n)}{\partial\lambda_{i}(n)\partial\lambda_{j}(n)}\biggr]\\
 &  & +\delta_{r(n-1)}\biggl[\frac{1}{2}G_{ki}^{+}(n)R_{j}^{kl}(n-1)T_{l}(n)+\frac{1}{2}G_{ki}^{+}(n)h^{kl}(n)G_{lj}^{-}(n-1)\\
 &  & \qquad+\frac{1}{2}R_{i}^{kl}(n)G_{kj}^{-}(n-1)T_{l}(n)\biggr]\\
 &  & +\delta_{r(n+1)}\biggl[\frac{1}{2}G_{ki}^{-}(n)R_{j}^{kl}(n+1)T_{l}(n+1)+\frac{1}{2}G_{ki}^{-}(n)h^{kl}(n+1)G_{lj}^{+}(n+1)\\
 &  & \qquad+\frac{1}{2}R_{i}^{kl}(n)G_{kj}^{+}(n+1)T_{l}(n+1)\biggr]
\end{eqnarray*}

Since $\Theta$, $\Psi$ and $g$ are known functions of $\lambda$,
the second order Taylor series is defined. One can now identify $G$,
$H^{+}$, $H^{-}$ and $K$ for initialising the Gaussian recursion
detailed above. Note that the $H^{\pm}$ are double the $\delta_{r(n\pm1)}$
coefficients above because when constructing (\ref{Marginal}) for
$k=m$ one needs to take into account equal contributions from $\frac{1}{\alpha}\frac{\partial^{2}S}{\partial\lambda_{i}(n)\partial\lambda_{j}(n+1)}$
and from $\frac{1}{\alpha}\frac{\partial^{2}S}{\partial\lambda_{j}(n+1)\partial\lambda_{i}(n)}$.
Explicitly we have
\begin{eqnarray*}
G_{ij}^{r}(m) & = & -V_{kij}(r)h^{kl}(r)T_{l}(r)+\frac{1}{2}G_{ki}^{+}(r)R_{j}^{kl}(r)T_{l}(r)\\
 &  & +\frac{1}{2}G_{ki}^{+}(r)h^{kl}(r)G_{lj}^{+}(r)+\frac{1}{2}U_{ij}^{kl}(r)T_{k}(r)T_{l}(r)+\frac{1}{2}R_{i}^{kl}(r)G_{kj}^{+}(r)T_{l}(r)\\
 &  & -V_{kij}(r)h^{kl}(r+1)T_{l}(r+1)+\frac{1}{2}G_{ki}^{-}(r)R_{j}^{kl}(n)T_{l}(r+1)\\
 &  & +\frac{1}{2}G_{ki}^{-}(r)h^{kl}(r+1)G_{lj}^{-}(r)+\frac{1}{2}U_{ij}^{kl}(r)T_{k}(r+1)T_{l}(r+1)\\
 &  & +\frac{1}{2}R_{i}^{kl}(r)G_{kj}^{-}(r)T_{l}(r+1)+\frac{\partial^{2}\Psi(r)}{\partial\lambda_{i}(r)\partial\lambda_{j}(r)}\\
H_{ij}^{-r}(m) & = & G_{ki}^{+}(r)R_{j}^{kl}(r-1)T_{l}(r)+G_{ki}^{+}(r)h^{kl}(r)G_{lj}^{-}(r-1)+R_{i}^{kl}(r)G_{kj}^{-}(r-1)T_{l}(r)\\
H_{ij}^{+r}(m) & = & G_{ki}^{-}(r)R_{j}^{kl}(r+1)T_{l}(r+1)+G_{ki}^{-}(r)h^{kl}(r+1)G_{lj}^{+}(r+1)+R_{i}^{kl}(r)G_{kj}^{+}(r+1)T_{l}(r+1)\\
K_{i}^{r}(m) & = & \frac{\partial\Psi(r)}{\partial\lambda_{i}(r)}+G_{ki}^{+}(r)h^{kl}(r)T_{l}(r)+\frac{1}{2}R_{i}^{kl}(r)T_{k}(r)T_{l}(r)\\
 &  & +G_{ki}^{-}(r)h^{kl}(r+1)T_{l}(r+1)+\frac{1}{2}R_{i}^{kl}(r)T_{k}(r+1)T_{l}(r+1)
\end{eqnarray*}

\subsection{Endpoint conditions}

At $n=1$ we need to define $T_{i}(1)$. This we do by assuming that
$\lambda_{i}(0)$ is fixed by some prescribed initial condition. At
$n=N$ we note that $T_{i}(N+1)$ does not occur in the discrete action
due to the truncation adopted. Thus all terms deriving from this drop
in the expressions above for the partial derivatives. Furthermore
in the second partial derivative $m=N+1$ is not present due to the
truncation so we obtain the following simpler expressions:
\begin{eqnarray*}
\frac{1}{\alpha}\frac{\partial S}{\partial\lambda_{i}(N)} & = & G_{ki}^{+}(N)h^{kl}(N)T_{l}(N)+\frac{1}{2}R_{i}^{kl}(N)T_{k}(N)T_{l}(N)+\frac{\partial\Psi(N)}{\partial\lambda_{i}(N)}
\end{eqnarray*}
\begin{eqnarray*}
\frac{1}{\alpha}\frac{\partial^{2}S}{\partial\lambda_{i}(N)\partial\lambda_{j}(m)} & = & \delta_{mN}\biggl[-V_{kij}(N)h^{kl}(N)T_{l}(N)+\frac{1}{2}G_{ki}^{+}(N)R_{j}^{kl}(N)T_{l}(N)+\frac{1}{2}G_{ki}^{+}(N)h^{kl}(N)G_{lj}^{+}(N)\\
 &  & +\frac{1}{2}U_{ij}^{kl}(N)T_{k}(N)T_{l}(N)+\frac{1}{2}R_{i}^{kl}(N)G_{kj}^{+}(N)T_{l}(N)+\frac{\partial^{2}\Psi(N)}{\partial\lambda_{i}(N)\partial\lambda_{j}(N)}\biggr]\\
 &  & +\delta_{m(N-1)}\biggl[\frac{1}{2}G_{ki}^{+}(N)R_{j}^{kl}(N-1)T_{l}(N)+\frac{1}{2}G_{ki}^{+}(N)h^{kl}(N)G_{lj}^{-}(N-1)\\
 &  & +\frac{1}{2}R_{i}^{kl}(N)G_{kj}^{-}(N-1)T_{l}(N)\biggr]
\end{eqnarray*}

\section{Endpoint considerations}

Path integrals calculated as above have fixed endpoints since the
action calculation detailed above demands this. For our application
however only the initial condition is fixed by the prescribed initial
trial density. Once a delta function for the initial consistency distribution
is prescribed then consistency distributions are determined for all
times since they satisfy a Euclidean Schrödinger equation. Thus the
final endpoint at $t=T$ must be sampled from a determined consistency
distribution. For a fixed final endpoint the consistency distribution
is proportional to the following sum of path weights
\begin{eqnarray}
\varrho(x_{0},x(T)) & = & {\displaystyle \sum_{\begin{array}{ccc}
\lambda(0) & = & x_{0}\\
\lambda(T) & = & x(T)
\end{array}}}\exp\left[-\Delta tS(\lambda)\right]\label{Boltzmann}
\end{eqnarray}

This equation is just a restatement of the generalized Boltzmann principle
for paths assumed in the present approach. By choosing a variety of
endpoints $x$ we can clearly calculate the relative magnitude of
the consistency distribution at time $T$ for various choices of $x$.
For many applications however $x$ is a vector and so only a sample
of possible such final time vectors will be available. The question
then becomes how to ensure that these follow the appropriate consistency
distribution. We proceed as follows: Construct a moderate sample $N$
of the $x(T)$ using a trial density\footnote{Such as an appropriate multivariate Gaussian calculated from a linearization
of the problem} $\phi(x)$. For each sample member perform the fixed endpoints MCMC
calculation as detailed previously and then calculate $\varrho$ using
(\ref{Boltzmann}). This calculation will be straightforward since
$S$ is central to the MCMC method. Consider now a small volume $dV$
surrounding $x$. The number of paths sampled will be 
\[
n_{t}(x)=N\phi(x)dV
\]

as opposed to the desired 
\[
n(x)=N\varrho(x_{0},x(T))dV
\]

Thus the contribution of each path ending in $x$ in calculating such
things as moments at any time along the path, will need to be reweighted
by the factor
\[
\Phi(x_{0},x)=\varrho(x_{0},x(T))/\phi(x)
\]

which may well vary significantly with each different sampled value
$x$ if the trial density differs from the actual density there.

\bibliographystyle{plain}
\bibliography{refs,bruce,/home/richard/Dropbox/Documents/neargaussian/biometrika2012/cv}

\end{document}